\documentstyle[11pt]{article}
\begin{document}

\title{Hypersurfaces with flat centroaffine metric and equations of
associativity}
\author{  {\Large Ferapontov E.V. } \\
    Department of Mathematical Sciences \\
    Loughborough University \\
    Loughborough, Leicestershire LE11 3TU \\
    United Kingdom \\
    e-mail: {\tt E.V.Ferapontov@lboro.ac.uk}
}
\date{}
\maketitle

\newtheorem{theorem}{Theorem}

\pagestyle{plain}

\maketitle

\begin{abstract}

It is demonstrated that hypersurfaces $M^n\subset A^{n+1}$ with a flat
centroaffine metric are governed
by a system of nonlinear PDEs known as the equations of associativity of
2-dimensional 
topological field theory. In the case of surfaces  $M^2\subset A^3$
this system reduces to a single
third order PDE,
$$
f_{xxx}f_{yyy}-f_{xxy}f_{xyy}=1,
$$
where $x$ and $y$ are the asymptotic coordinates on $M^2$.

\bigskip

Keywords: centroaffine geometry, equations of associativity.

\bigskip

Mathematical Subject Classification: 53A40, 53A20, 35Q58, 35L60.
\end{abstract}

\newpage
\section{Introduction}

Centroaffine geometry studies the properties of hypersurfaces
$M^n\subset A^{n+1}$
which are invariant under the action of the general linear group $GL(n+1)$.
Since  this action preserves the origin (translations are not allowed), one
can 
choose the position vector ${\bf r}$ of a hypersurface  as its
centroaffine normal. With this normalization, one readily introduces the
main 
centroaffine invariants: the
centroaffine metric, the centroaffine cubic form and the Chebyshev
covector  (see sect. 1; in our presentation we follow the monograph
\cite{Shirokov}, see also \cite{Scharlach}). These invariants satisfy the
compatibility 
conditions known as the
centroaffine Gauss-Codazzi equations. It is well-known that the
vanishing of the centroaffine cubic form characterises central
hyperquadrics, 
while the vanishing of the Chebyshev covector is characteristic for
proper affine hyperspheres. In this article we consider the third natural
class
of centroaffine hypersurfaces characterized by the flatness of the
centroaffine metric.

The main observation is that, in the flat coordinates
of the centroaffine metric,
the corresponding Gauss-Codazzi equations reduce to a system of nonlinear
PDEs 
known as the 
equations of associativity of two-dimensional topological field theory
\cite{Dubrovin}
(to be precise, WDVV equations without the quasihomogeneity condition). This
is explained in 
sect. 3. For surfaces  $M^2\subset A^3$  these equations reduce to a single
PDE of the Monge-Ampere type,
$$
f_{xxx}f_{yyy}-f_{xxy}f_{xyy}=1,
$$
where $x$ and $y$ are the asymptotic coordinates on $M^2$.

The geometry of surfaces $M^2\subset A^3$ with a flat centroaffine metric
is discussed in sect. 4 and 5 where, in particular, we establish some
geometric 
properties of characteristics of the above PDE.

The review  \cite{Dubrovin} provides an abundance of exact solutions of the
associativity equations. Some of them are listed in sect. 5 and 6.
The centroaffine geometry of the corresponding surfaces is currently under
investigation.

In sect. 7  a simple loop-group formulation of the linear system
governing surfaces with a flat centroaffine metric is proposed
which is due to A. Bobenko.

Since the centroaffine metric is proportional to the second fundamental
form, 
hypersurfaces with a flat centroaffine metric (and projective transforms
thereof) 
form a subclass of projective hypersurfaces with conformally flat second
fundamental form which was investigated in detail in \cite{Ak2}. Another
interesting class of projective surfaces $M^2\subset P^3$
(containing projective transforms of surfaces with a flat centroaffine
metric as a 
proper subclass), is characterised by the existence of a one-parameter
family of
asymptotic deformations
which induce constant rescalings of the projective metric and the Darboux
cubic form. 
These  generalizations are discussed in sect. 8.

\section{Hypersurfaces in centroaffine geometry}

The position vector ${\bf r}=(r^0, r^1, ..., r^n)$ of a hypersurface $M^n$
in centroaffine geometry
satisfies the linear system
\begin{equation}
{\bf r}_{ij}=\tilde {\Gamma} ^k_{ij} \ {\bf r}_k +g_{ij} \ {\bf r}, ~~~ i,
j=1, ..., n,
\label{centro}
\end{equation}
which is required to be compatible of the rank $n+1$
(here ${\bf r}_i=\partial {\bf r}/
\partial u^i, \ {\bf r}_{ij}=\partial {\bf r}/
\partial u^i \partial u^j$ where $u^1, ..., u^n$ are local coordinates on
$M^n$).
 Under reparametrizations of $M^n$, the objects $g_{ij}$ and
$\tilde {\Gamma} ^k_{ij}$ transform as the components of a pseudo-Riemannian
metric
(assumed nondegenerate) and Christoffel's symbols of an affine connection
(called the affine connection of the first kind;
we emphasize that it
is {\it not} the Levi-Civita connection of $g_{ij}$). The conformal class of
$g_{ij}$ is nothing but the second fundamental form of $M^n$. The
compatibility conditions of
system (\ref{centro}) are of the form \cite{Shirokov}
\begin{equation}
\begin{array}{c}
\tilde \nabla _k g_{ij}=\tilde \nabla _j g_{ik} \\
\ \\
\tilde R^s_{ijk}=g_{ik}\delta^s_j-g_{ij}\delta^s_k
\end{array}
\label{comp1}
\end{equation}
where $\tilde \nabla$ denotes covariant differentiation in the connection
$\tilde \Gamma$, and $\tilde R^s_{ijk}=
\partial_k \tilde \Gamma ^s_{ij}-\partial_j \tilde \Gamma ^s_{ik}
+\tilde \Gamma^p_{ij}\tilde \Gamma^s_{pk}-
\tilde \Gamma^p_{ik}\tilde \Gamma^s_{pj}$ is the curvature tensor. For our
purposes, it will
be more convenient to rewrite the compatibility conditions (\ref{comp1}) in
terms of the Levi-Civita connection $\Gamma ^k_{ij}$ of the metric $g_{ij}$
and the
difference $(1, 2)$-tensor $h^k_{ij}$ defined by
$$
\tilde \Gamma ^k_{ij}=\Gamma ^k_{ij}+h^k_{ij}.
$$
In this notation, equations (\ref{comp1}) can
be cast into the form \cite{Shirokov}
\begin{equation}
\begin{array}{c}
h^s_{ij} g_{sk}=h^s_{ik} g_{sj}, ~~~
\nabla _k h^s_{ij}= \nabla _j h^s_{ik}, \\
\ \\
R^s_{ijk}+h^p_{ij}h^s_{pk}-h^p_{ik}h^s_{pj}
=g_{ik}\delta^s_j-g_{ij}\delta^s_k
\end{array}
\label{comp2}
\end{equation}
where $\nabla$ denotes covariant differentiation in the Levi-Civita
connection 
$\Gamma$, and  $R^s_{ijk}$ is the curvature tensor. Equation
$(\ref{comp2})_1$ implies that the tensor $h_{ijk}=h^s_{ij} g_{sk}$
is totally symmetric, defining the centroaffine cubic form
of the hypersurface $M^n$. The centroaffine metric
$$
M=g_{ij}\ du^idu^j
$$ 
and the centroaffine cubic form
$$
C=h_{ijk}\ du^idu^jdu^k
$$ 
satisfying the compatibility conditions (\ref{comp2})
characterize a hypersurface uniquely. Another useful invariant of a
hypersurface in centroaffine geometry is the so-called
Chebyshev covector,
$$
T=T_i\ du^i=\frac{1}{n}h^s_{is}\ du^i,
$$
which, as readily follows from
$(\ref{comp2})_2$, is always a gradient: $dT=0$. It is known that the
vanishing of 
the cubic form characterizes central hyperquadrics \cite{Shirokov}.
The proper affine hyperspheres can be defined as hypersurfaces with the
vanishing Chebyshev covector. The main object of our study is the third
natural class of centroaffine
hypersurfaces, namely, hypersurfaces with a flat centroaffine metric.
Notice that the centroaffine metric
can be equivalently defined by the formula
$$
M=\frac{\vert d^2{\bf r}\wedge {\bf r}_1\wedge ... \wedge {\bf r}_n\vert}
{\vert {\bf r}\wedge {\bf r}_1\wedge ... \wedge {\bf r}_n\vert}.
$$
For hypersurfaces defined explicitly as graphs of functions,
$u^{n+1}=f(u^1, ..., u^n)$, the centroaffine metric has the form
$$
M=-\frac{d^2f}{F}=-\frac{d^2F}{F},
$$ 
where $F=f_iu^i-f$ is the Legendre transform of $f$, and $d^2F$
denotes the second differential. Some global aspects of
Riemannian metrics  of the form $-d^2F/F$ were discussed in \cite{Loftin}
and 
\cite{Loewner}. The condition of flatness of the metric $-d^2F/F$ results in
a  system of third order nonlinear PDEs for the function $F$ (which looks
quite
formidable even in the simplest case $n=2$). As we demonstrate below,
this system is nothing but the equations of associativity (WDVV equations).

{\bf Remark. } System (\ref{centro}) appears naturally in the theory of
multi-dimensional 
Schwarzian derivatives. In this setting, one is interested in the mapping
from $M^n$ into the projective space $P^n$ defined as
$(r^0 : r^1 : ... : r^n)$ in
homogeneous coordinates of $P^n$. We refer to \cite{Yoshida} for the
details.
In centroaffine geometry, non-constant
rescalings of the position vector ${\bf r}$ are not allowed: they
result in non-equivalent hypersurfaces.

\section{Hypersurfaces with flat centroaffine metric
and associativity equations}

For hypersurfaces with a flat centroaffine metric  one can
choose the parametrization such that $g_{ij}=\eta_{ij}=const$, so that
$\Gamma =0,
\ \nabla _k=\partial _k$  and equations (\ref{comp2}) take the form
\begin{equation}
\begin{array}{c}
h^s_{ij} \eta_{sk}=h^s_{ik} \eta_{sj}, ~~~
\partial _k h^s_{ij}= \partial _j h^s_{ik}, \\
\ \\
h^p_{ij}h^s_{pk}-h^p_{ik}h^s_{pj}
=\eta _{ik}\delta^s_j-\eta _{ij}\delta^s_k.
\end{array}
\label{comp3}
\end{equation}
The first two equations  imply the existence of the potential $f$ such that
$$
h^s_{ij}=\eta^{sk} \partial_i\partial_j\partial_k f
$$
resulting in the following formulae for the centroaffine invariants:

\noindent centroaffine metric
$$
M=\eta_{ij}\ du^idu^j;
$$

\noindent centroaffine cubic form
$$
C=d^3f;
$$

\noindent Chebyshev covector
$$
T=\frac{1}{n}\ d(\triangle f)
$$
where $\triangle$ denotes the Laplacian in the metric $\eta_{ij}$.

\noindent In terms of $f$,  the third set of equations (\ref{comp3}) results
in a nonlinear system of PDEs
known as the equations of associativity in two-dimensional topological field
theory
\cite{Dubrovin}. To see this,  we introduce one extra coordinate $u^0$ and
the function
$$
F(u^0, u^1, ..., u^n)=\frac{1}{6} (u^0)^3+\eta_{ij}u^0u^iu^j+f(u^1, ...,
u^n).
$$
One can readily verify the following properties, which are in fact the
axioms of
the associativity equations (in what follows $\alpha, \beta, \gamma
 =0, 1, ..., n ~ {\rm and} ~ i, j, k =1, ..., n$).

\bigskip

\noindent 1. The matrix $\eta_{\alpha \beta}= \frac{\partial^3 F}
{\partial_0\partial_{\alpha} \partial_{\beta}}$ is constant and
nondegenerate,
indeed, 
$$
\eta_{\alpha \beta}=
\left(
\begin{array}{cc}
1&0 \\
\ \\
0&\eta_{ij}
\end{array} \right).
$$

\medskip

\noindent 2. The objects $c^{\alpha}_{\beta \gamma}=
\eta^{\alpha \delta}\frac{\partial^3F}{\partial_ {\beta} \partial_ {\gamma}
 \partial_ {\delta}}
$ are the structure constants of a commutative associative algebra. Indeed,
the
only nonzero among them are
$$
c^0_{00}=1, ~~~ c^0_{ij}=\eta_{ij}, ~~~ c^i_{0j}=c^i_{j0}=\delta^i_j, ~~~
c^i_{jk}=h^i_{jk},
$$
and equations (\ref{comp3}) imply the associativity.

The function $F$ defines a Frobenius structure \cite{Dubrovin}
on the manifold with coordinates
$u^0, u^1, ..., u^n$, the main ingredients of which are the flat metric
$$
(du^0)^2+\eta_{ij}\ du^idu^j
$$
and the symmetric 3-tensor
$$
d^3F=(du^0)^3+\eta_{ij} du^0du^idu^j+d^3f.
$$
Notice that the centroaffine metric $M$ and the cubic form $C$ can be
obtained 
by resticting these objects to a hyperplane $u^0=const$ in the spirit of
\cite{Strachan}. Thus, hypersurfaces with flat centroaffine metric
carry an intrinsic induced Frobenius structure.
Moreover,  one can incorporate the
spectral parameter $\lambda$ into the equations for the position vector,
\begin{equation}
{\bf r}_{ij}=\lambda h^k_{ij} \ {\bf r}_k +\lambda^2 \eta_{ij} \ {\bf r},
\label{Lax}
\end{equation}
without violating the compatibility conditions (\ref{comp3}).
In arbitrary parametrization, the linear system (\ref{Lax}) takes the form
\begin{equation}
{\bf r}_{ij}=(\Gamma^k_{ij}+\lambda h^k_{ij}) \ {\bf r}_k +
\lambda^2 g_{ij} \ {\bf r},
\label{Lax0}
\end{equation}
where $\Gamma$ is the Levi-Civita connection of the flat metric $g$.
System (\ref{Lax0}) defines a one-parameter family of hypersurfaces
$M^n_{\lambda}$  with  centroaffine metrics
$$
M=\lambda^2 g_{ij}\ du^idu^j
$$
and  centroaffine cubic forms
$$
C= \lambda^{3} h_{ijk}\ du^idu^jdu^k.
$$
Therefore, hypersurfaces with a flat centroaffine metric can be equivalently
characterized
as hypersurfaces possessing nontrivial deformations which induce rescalings
of the 
centroaffine metric and the centroaffine cubic form. For $\lambda=1$ we
recover the original hypersurface $M^n$, the case
 $\lambda=0$ corresponds to a hyperplane, and $\lambda =-1$ results in the
dual hypersurface.

{\bf Remark.} One can consider the linear system (\ref{Lax0}) without
imposing 
the restriction that $\Gamma$ is the Levi-Civita connection of $g$. In this
case the 
compatibility conditions imply
\begin{equation}
\begin{array}{c}
\nabla_k g_{ij}=\nabla_j g_{ik}, ~~~
\nabla_k h^p_{ij}=\nabla_j h^p_{ik}, ~~~
h^s_{ij}g_{sk}=h^s_{ik}g_{sj}, \\
\ \\
h^p_{ij}h^s_{pk}-h^p_{ik}h^s_{pj}
=g_{ik}\delta^s_j-g_{ij}\delta^s_k
\end{array}
\label{gen}
\end{equation}
where $\nabla$ denotes covariant differentiation in the connection $\Gamma$.
Moreover, the connection $\Gamma$ must be flat. This
class of surfaces will be investigated in some more detail elsewhere.

{\bf Examples.} In the case $n=2$ one can either choose  $\eta = 2\ dxdy$
(hyperbolic surfaces) or $\eta = dx^2+dy^2$ (convex surfaces).
Here $x=u^1$ and $y=u^2$ are asymptoticcoordinates on a surface $M^2$.
In the hyperbolic case,
equations (\ref{Lax})  take the form
\begin{equation}
\begin{array}{c}
{\bf r}_{xx}=\lambda f_{xxx}{\bf r}_y+\lambda f_{xxy}{\bf r}_x, \\
\ \\
{\bf r}_{xy}=\lambda f_{xxy}{\bf r}_y+
\lambda f_{xyy}{\bf r}_x+\lambda^2 {\bf r}, \\
\ \\
{\bf r}_{yy}=\lambda f_{xyy}{\bf r}_y+\lambda f_{yyy}{\bf r}_x,
\end{array}
\label{Lax1}
\end{equation}
with the compatibility condition
\begin{equation}
f_{xxx}f_{yyy}-f_{xxy}f_{yyx}=1.
\label{ass1}
\end{equation}
In a similar form  (however, without a spectral parameter) the linear
problem
  (\ref{Lax1}) was presented in
\cite{Dryuma}. Introducing
$\psi=(\lambda {\bf r}, {\bf r}_x, {\bf r}_y)^T$,
one can readily rewrite (\ref{Lax1}) in a matrix form
\begin{equation}
\psi_x=\lambda \left(
\begin{array}{ccc}
0&1&0\\
0&f_{xxy}&f_{xxx}\\
1&f_{xyy}&f_{xxy}
\end{array} \right) \psi, ~~~~
\psi_y=\lambda \left(
\begin{array}{ccc}
0&0&1\\
1&f_{xyy}&f_{xxy}\\
0&f_{yyy}&f_{xyy}
\end{array} \right) \psi
\label{trace}
\end{equation}
which coincides with the one from \cite{Dubrovin}.
Similarly, in the convex case we have
\begin{equation}
\begin{array}{c}
{\bf r}_{xx}=\lambda f_{xxx}{\bf r}_x+\lambda f_{xxy}{\bf r}_y+
\lambda^2{\bf r}, \\
\ \\
{\bf r}_{xy}=\lambda f_{xxy}{\bf r}_x+
\lambda f_{xyy}{\bf r}_y, \\
\ \\
{\bf r}_{yy}=\lambda f_{xyy}{\bf r}_x+\lambda f_{yyy}{\bf r}_y+\lambda^2
{\bf r},
\end{array}
\label{Lax2}
\end{equation}
with the compatibility condition
\begin{equation}
f_{xxy}^2+f_{xyy}^2-f_{xxx}f_{xyy}-f_{yyy}f_{xxy}=1.
\label{ass2}
\end{equation}
Particular examples of surfaces with a flat centroaffine metric were
discussed in 
\cite{Scharlach} (surfaces of revolution with the centroaffine metric of
constant curvature) and
\cite{Wang} (centroaffinely flat Chebyshev surfaces).

Introducing  the variables $f_{xxx}=p, \ f_{xxy}=a, \ f_{xyy}=b$ and
$\ f_{yyy}=q$,
\cite{Mokhov}, \cite{Fer96}, one can rewrite (\ref{ass1}) in the form
\begin{equation}
p_y=a_x, ~~~ a_y=b_x,  ~~~ b_y=q_x, ~~~ pq-ab=1
\label{hydro0}
\end{equation}
which is equivalent to the quasilinear system
\begin{equation}
\left( \begin{array}{c}
p\\
a\\
b 
\end{array}\right)_y=U\ \left( \begin{array}{c}
p\\
a\\
b 
\end{array} \right)_x, ~~~~~ U=\left( \begin{array}{ccc}
0&1&0\\
0&0&1\\
-\frac{1+ab}{p^2}&\frac{b}{p}&\frac{a}{p}
\end{array}
\right)
\label{hydro}
\end{equation}
investigated in detail in \cite{Fer96}. The characteristics of this system
are defined by  $dx+\lambda \ dy=0$
where $\lambda$ is a root of the
characteristic polynomial of the matrix $U$,
$$
p\lambda^3-a\lambda^2-b\lambda+q=0.
$$
Equivalently, characteristics can be defined as the null-curves of the cubic
form
$$
p\ dx^3+a \ dx^2dy-b\ dx dy^2 - q\ dy^3,
$$
indeed, $\lambda=-dx/dy$. The geometry of this cubic form is
clarified in section 4.

\section{Surfaces in centroaffine 3-space}

In asymptotic coordinates $x$ and $y$
(for definiteness, we consider hyperbolic case only),
 equations  (\ref{centro}) take the form
\begin{equation}
\begin{array}{c}
{\bf r}_{xx}=(a+u_x)\ {\bf r}_x+p\ {\bf r}_y, \\
\ \\
{\bf r}_{xy}=b\ {\bf r}_x+a\ {\bf r}_y+e^{u} \ {\bf r}, \\
\ \\
{\bf r}_{yy}=q\ {\bf r}_x+(b+u_y)\ {\bf r}_y,
\end{array}
\label{centro2}
\end{equation}
with the compatibility conditions
\begin{equation}
\begin{array}{c}
a_y=b_x, ~~~ a_x-p_y=au_x+pu_y, ~~~ b_y-q_x=qu_x+bu_y, \\
\ \\
pq-ab=e^u(K+1),
\end{array}
\label{comp4}
\end{equation}
where $K=-e^{-u}u_{xy}$ is the Gaussian curvature of the corresponding
centroaffine metric
$$
M=2e^u\ dxdy.
$$
Notice that in the flat case ($K=0, \ u=0$)
equations (\ref{comp4}) coincide with (\ref{hydro0}).
Since the only nonzero Christoffel symbols of the centroaffine metric $M$
are
$\Gamma ^1_{11}=u_x$ and $\Gamma ^2_{22}=u_y$, we have
$h^1_{11}=a, \ h^2_{11}=p, \ h^1_{12}=b, \ h^2_{12}=a, \ h^1_{22}=q$ and
$ h^2_{22}=b$, so that the centroaffine cubic form $C$ and the Chebyshev
covector
$T$ take the forms
$$
C=e^u(p\ dx^3+3a\ dx^2dy+3b\ dxdy^2+q\ dy^3)
$$
and
$$
T=a\ dx+b\ dy,
$$
respectively. 

{\bf Remark}. We point out that surfaces with the zero Chebyshev covector
are the 
proper affine spheres, indeed, the substitution of  $a=b=0$ into
$(\ref{comp4})_2$ and  $(\ref{comp4})_3$ implies $p=q=e^{-u}$, so that
$(\ref{comp4})_4$ results in the Tzitzeica equation $u_{xy}=e^u-e^{-2u}$.

Let us consider the form
$$
C-2\ TM=e^u(p\ dx^3-a\ dx^2dy-b\ dxdy^2+q\ dy^3)
$$
which, at a generic point of our surface, has three distinct zero
directions. 
The three directions conjugate to them are defined by the equation
\begin{equation}
p\ dx^3+a\ dx^2dy-b\ dxdy^2-q\ dy^3=0,
\label{char}
\end{equation}
(recall that two directions are called conjugate if they are
orthogonal with respect to the metric $M$). The null-curves of
(\ref{char}) define a foliation of our surface by three families of curves,
which we will call the {\it characteristic 3-web} for the reason which
is clarified at the end of section 3.

\begin{theorem}

For surfaces with a flat centroaffine metric, the characteristic
3-web is hexagonal.

\end{theorem}

\centerline{\bf Proof:}

A direct calculation of the curvature form $\Omega$
of the characteristic 3-web results in the formula $\Omega=d\omega$, where
the connection form $\omega$ is
$$
\omega=-\frac{3e^u}{r^2-4st}
\left( (rK_x-2tK_y)\ dx+(rK_y-2sK_x)\ dy\right),
$$
here $r=ab-9pq, \ s= b^2+3aq, \ t=a^2+3bp$. If $K=0$, the
curvature form $\Omega$ is zero, implying the hexagonality of the
characteristic 3-web.

{\bf Remark.} The class of surfaces with a hexagonal characteristic 3-web
includes surfaces with a flat centroaffine metric as a proper subclass.
For instance, the form of $\omega$ implies that  surfaces with the
centroaffine metric of constant curvature ($K$=const)  also have a
hexagonal characteristic web. In would be of interest to discuss this class
of surfaces in some more detail. The simplest examples of surfaces with
centroaffine metric of constant curvature
$K=-1$ are central quadrics
(which correspond to $a=b=p=q=0$ in the equations (\ref{comp4})). Since in
this case
(\ref{char}) is identically zero,
the characteristic 3-web is indeterminate.

\section{Other parametrizations}

In asymptotic coordinates $x$ and $y$, the position vector $\bf r$ of a
surface with flat centroaffine metric satisfies the equations
\begin{equation}
\begin{array}{c}
{\bf r}_{xx}=\lambda a\ {\bf r}_x+\lambda p\ {\bf r}_y, \\
\ \\
{\bf r}_{xy}=\lambda b\ {\bf r}_x+\lambda a\ {\bf r}_y+\lambda^2 \ {\bf r},
\\
\ \\
{\bf r}_{yy}=\lambda q\ {\bf r}_x+\lambda b\ {\bf r}_y,
\end{array}
\label{paramass}
\end{equation}
the compatibility conditions of which coincide with (\ref{hydro0}).
After the reparametrization $(x, y)\to (t, y)$ where $t$ is defined by
\begin{equation}
dt=pdx+ady,
\label{t}
\end{equation}
equations (\ref{paramass}) take the form
\begin{equation}
\begin{array}{c}
{\bf r}_{tt}= (\frac{C_t}{C}-2\lambda B)\ {\bf r}_t     +\lambda C\ {\bf
r}_y\\
\ \\
{\bf r}_{ty}=(\frac{B_t}{C}-\lambda A)\ {\bf r}_t+ \lambda^2 C \ {\bf r}, \\
\ \\
{\bf r}_{yy}=(\lambda+\frac{B_y}{C})\ {\bf r}_t-\lambda A\ {\bf r}_y+
2\lambda^2 B\ {\bf r},
\end{array}
\label{paramass1}
\end{equation}
where we have introduced the notation
\begin{equation}
A={a^2}/{p}-b, ~~~ B=-{a}/{p}, ~~~ C={1}/{p}.
\label{connection}
\end{equation}
The compatibility conditions of  equations (\ref{paramass1}) take the form
$$
A_t=B_y, ~~~ B_t=C_y, ~~~ C_t=(B^2-AC)_y.
$$
Written in terms of the potential $F$ such that $A=F_{yyy}, ~ B=F_{tyy}$ and
$C=F_{tty}$, this system reduces to a single PDE
\begin{equation}
F_{ttt}=F_{tyy}^2-F_{tty}F_{yyy},
\label{ass3}
\end{equation}
which is yet another form of the associativity equations \cite{Dubrovin}.
Thus, both equations (\ref{ass1}) and (\ref{ass3}) describe one in the
same class of surfaces, however, in different parametrizations. Notice that,
in terms of the potentials $f(x, y)$ and $F(x, t)$, the formulae (\ref{t})
and 
(\ref{connection}) take the form
$$
t=f_{xx}, ~~~ F_{yyy}=\frac{f_{xxy}^2}{f_{xxx}}-f_{xyy}, ~~~
F_{tyy}=-\frac{f_{xxy}}{f_{xxx}}, ~~~ F_{tty}=\frac{1}{f_{xxx}},
~~~ F_{ttt}=\frac{f_{xyy}}{f_{xxx}},
$$
and, after being integrated once, simplify to
\begin{equation}
f_{xx}=t, ~~~ f_{xy}=-F_{yy}, ~~~ f_{yy}=F_{tt}, ~~~ x=F_{ty}.
\label{trans}
\end{equation}
This
transformation between (\ref{ass1}) and (\ref{ass3}) was first proposed in
\cite{Fer96}. Taking, for instance, the known polynomial solution
of the equation (\ref{ass3}),
$$
F(t, y)=\frac{t^2y^2}{4}+\frac{t^5}{60},
$$
 we readily arrive at the formulae
$$
f_{xx}=t, ~~~ f_{xy}=-\frac{t^2}{2}, ~~~ f_{yy}=\frac{y^2}{2}+\frac{t^3}{3},
~~~ 
x=ty,
$$
which, after the substitution of $t=x/y$ into the first three equations,
give a solution of  (\ref{ass1}):
$$
f(x, y)=\frac{x^3}{6y}+\frac{y^4}{24}.
$$
Here we list some further solutions of the equation
$f_{xxx}f_{yyy}-f_{xxy}f_{yyx}=1$ which can be obtained by applying the
transformation
(\ref{trans}) to the known solutions of (\ref{ass3})
as found in \cite{Dubrovin}:

\vspace{3ex}
\begin{tabular}{|l|l|}                                         \hline
                                          &                     \\ [-2ex]
 $F_{ttt}=F_{tyy}^2-F_{tty}F_{yyy}$ &  $f_{xxx}f_{yyy}-f_{xxy}f_{yyx}=1$ \\
                                                                 [1ex]
\hline
                                          &                      \\ [-2ex]
$F=\frac{t^2y^2}{4}+\frac{t^5}{60}$       &
$f=\frac{x^3}{6y}+\frac{y^4}{24}$ \\
                                                                 [2ex]
\hline
                                          &                      \\ [-2ex]
$F=ye^t-\frac{y^4}{24}$
& $f=\frac{xy^3}{6}+\frac{x^2 \ln x}{2}-\frac{3x^2}{4}$ \\
                                                                 [2ex]
\hline
                                          &                      \\ [-2ex]
$F=\frac{y^2e^t}{4}+\frac{e^{2t}}{32}-\frac{y^4}{48}$
&  $f=\frac{xy^3}{12}-\frac{x^2\ln y}{2}+\frac{x^2 \ln x}{2}-\frac{3x^2}{4}$
\\
                                                                 [2ex]
\hline
                                          &                      \\ [-2ex]
$F=\frac{t^2\ln y}{2}$
& $f=\frac{yx^3}{6}+\frac{y^2 \ln y}{2}-\frac{3y^2}{4}$    \\
                                                                 [2ex]
\hline
                                          &                      \\ [-2ex]
$F=\frac{y^3t}{6}+\frac{y^2t^3}{6}+\frac{t^7}{210}$
& $f=\frac{4y^2}{15}\left( \frac{x}{y}-\frac{y}{2}\right)^{5/2}$ \\
                                                                 [2ex]
\hline
                                          &                      \\ [-2ex]

$F=\frac{y^3t^2}{6}+\frac{y^2t^5}{20}+\frac{t^{11}}{3960}$
& $f=\frac{y^4t^3}{6}+\frac{7y^3t^6}{30}+\frac{4y^2t^9}{45}+\frac{y^5}{60}$,
\ where \ 
$x=y^2t+\frac{yt^4}{2}$ \\
                                                                   [1ex]
\hline
\end{tabular}

\vspace{3ex}

\noindent Notice that the second and fourth solutions in the right column
differ by the interchange of $x$ and $y$,
which is an obvious symmetry of the corresponding PDE. The last solution is
implicit.
The geometry of the corresponding surfaces is currently under investigation.

\bigskip

\noindent {\bf The characteristic parametrization}. Introducing the
variables
$w^1, w^2, w^3$ by the formulae
$$ 
w^1+w^2+w^3=2a, ~~~ w^1w^2+w^1w^3+w^2w^3=a^2-pb, ~~~ w^1w^2w^3=p,
$$
one can rewrite equations (\ref{hydro0}) in a symmetric form \cite{Fer96},
\begin{equation}
w^1_y=\left(\frac{w^1-w^2-w^3}{2w^2w^3}\right)_x, ~~~
w^2_y=\left(\frac{w^2-w^1-w^3}{2w^1w^3}\right)_x, ~~~
w^3_y=\left(\frac{w^3-w^1-w^2}{2w^1w^2}\right)_x.
\label{sym}
\end{equation}
This system  possesses three characteristic conservation laws
$$
\begin{array}{c}
(w^2-w^3)\left(dx+\frac{w^1+w^2-w^3}{2w^1w^2w^3} dy\right), ~~~~~
(w^1-w^3)\left(dx+\frac{w^1+w^3-w^2}{2w^1w^2w^3} dy\right), \\
\ \\
(w^2-w^1)\left(dx+\frac{w^2+w^3-w^1}{2w^1w^2w^3} dy\right), ~~~
\end{array}
$$
(we call a conservation law `characteristic' if its null curves are
characteristics of the system (\ref{sym})). After introducing the
characteristic
coordinates $\xi$ and $\eta$,
$$
\begin{array}{c}
d\xi=(w^2-w^3)\left(dx+\frac{w^1+w^2-w^3}{2w^1w^2w^3} dy\right), ~~~~~
d\eta =(w^1-w^3)\left(dx+\frac{w^1+w^3-w^2}{2w^1w^2w^3} dy\right), \\
\ \\
d\xi-d\eta=(w^2-w^1)\left(dx+\frac{w^2+w^3-w^1}{2w^1w^2w^3} dy\right),
\end{array}
$$
equations (\ref{paramass}) for the position vector take the form
\begin{equation}
\begin{array}{c}
{\bf r}_{\xi \xi}=(\Gamma^1_{11}+\lambda h^1_{11})\ {\bf r}_{\xi}
+(\Gamma^2_{11}+\lambda h^2_{11})\ {\bf r}_{\eta}+\lambda^2 E\ {\bf r}, \\
\ \\
{\bf r}_{\xi \eta}=(\Gamma^1_{12}+\lambda h^1_{12})\ {\bf r}_{\xi}
+(\Gamma^2_{12}+\lambda h^2_{12})\ {\bf r}_{\eta}+\lambda^2 F\ {\bf r}, \\
\ \\
{\bf r}_{\eta \eta}=(\Gamma^1_{22}+\lambda h^1_{22})\ {\bf r}_{\xi}
+(\Gamma^2_{22}+\lambda h^2_{22})\ {\bf r}_{\eta}+\lambda^2 G\ {\bf r},
\end{array}
\label{charcoord}
\end{equation}
where the coefficients of the centroaffine metric $M=Ed\xi^2+2F d\xi d\eta
+Gd\eta^2$
are given by
$$
\begin{array}{c}
E=\frac{w^1w^2w^3(w^2-w^1-w^3)}{(w^2-w^1)^2(w^2-w^3)^2}, ~~~~~
F=\frac{w^1w^2(w^3)^2}{(w^1-w^2)^2(w^1-w^3)(w^2-w^3)}, \\
\ \\
G=\frac{w^1w^2w^3(w^1-w^2-w^3)}{(w^1-w^2)^2(w^1-w^3)^2}, ~~~
\end{array}
$$
and $\Gamma^i_{jk}$ are the components of the corresponding Levi-Civita
connection.
The coefficients $h^i_{jk}$ are of the form
$$
\begin{array}{c}
h^1_{11}=\frac{(w^2)^2-w^1w^2-w^1w^3-w^2w^3}{(w^2-w^1)(w^2-w^3)}, ~~~~
h^2_{22}=\frac{(w^1)^2-w^1w^2-w^1w^3-w^2w^3}{(w^1-w^2)(w^1-w^3)}, \\
\ \\
h^1_{12}=-\frac{w^2w^3}{(w^1-w^2)(w^1-w^3)}, ~~~~ h^2_{11}=h^1_{22}=0, ~~~~
h^2_{12}=-\frac{w^1w^3}{(w^2-w^1)(w^2-w^3)}.
\end{array}
$$
We point out that the components $E, F$ and $G$ of the flat  metric $M$
satisfy the
identity
$$
(EG-F^2)^2+F(E+F)(G+F)=0.
$$
In  coordinates $\xi $ and $\eta$,
the characteristic 3-web is defined by the equation $d\xi d\eta (d\xi
-d\eta)=0$,
and is manifestly hexagonal.
The compatibility conditions of equations (\ref{charcoord}) take the form
$$
\begin{array}{c}
\frac{w^3-w^2}{w^1}w^1_{\xi}=\frac{w^3-w^1}{w^2}w^2_{\eta}, \\
\ \\
\frac{w^2-w^3}{w^1}(w^1_{\xi}+w^1_{\eta})=\frac{w^1-w^2}{w^3}w^3_{\eta}, \\
\ \\
\frac{w^1-w^3}{w^2}(w^2_{\xi}+w^2_{\eta})=\frac{w^2-w^1}{w^3}w^3_{\xi},
\end{array}
$$
and, upon the introduction of variables $u^i=1/w^i$,  simplify to
\begin{equation}
\begin{array}{c}
(u^3-u^2)u^1_{\xi}=(u^3-u^1)u^2_{\eta}, \\
\ \\
(u^2-u^3)(u^1_{\xi}+u^1_{\eta})=(u^1-u^2)u^3_{\eta}, \\
\ \\
(u^1-u^3)(u^2_{\xi}+u^2_{\eta})=(u^2-u^1)u^3_{\xi}.
\label{26}
\end{array}
\end{equation}
System (\ref{26}) can be represented in the equivalent exterior form
$$
\omega^1\wedge d\xi=0, ~~~ \omega^2\wedge d\eta=0, ~~~
\omega^3\wedge (d\xi-d\eta)=0,
$$
with the 1-forms
$$
\begin{array}{c}
\omega^1=\frac{(u^2-u^3)du^1+(u^1-u^3)du^2+
(u^2-u^1)du^3}{2(u^2-u^3)\sqrt {(u^2-u^1)(u^3-u^1)}}, \\
\ \\
\omega^2=\frac{(u^2-u^3)du^1+(u^1-u^3)du^2+
(u^1-u^2)du^3}{2(u^3-u^1)\sqrt {(u^2-u^1)(u^2-u^3)}}, \\
\ \\
\omega^3=\frac{(u^2-u^3)du^1+(u^3-u^1)du^2+
(u^2-u^1)du^3}{2(u^2-u^1)\sqrt {(u^3-u^1)(u^2-u^3)}} \\
\end{array}
$$
satisfying the structure equations of the $SO(2, 1)$ group,
\begin{equation}
d\omega^1=\omega^2\wedge \omega^3, ~~~
d\omega^2=\omega^3\wedge \omega^1, ~~~
d\omega^3=\omega^2\wedge \omega^1.
\label{structure}
\end{equation}
Therefore, one can set
$$
\omega^1=p^1d\xi, ~~~ \omega^2=p^2d \eta, ~~~ \omega^3=p^3(d\xi-d\eta),
$$
and the substitution into (\ref{structure}) implies the 3-wave system
\begin{equation}
p^1_{\eta}=p^2p^3, ~~~ p^2_{\xi}=p^1p^3, ~~~ p^3_{\xi}+p^3_{\eta}=p^1p^2.
\label{wave}
\end{equation}
It should be pointed out that the equivalence of the associativity
equations and the N-wave system was first observed by Dubrovin
\cite{Dubrovin}.
The explicit sequence of transformations shown above (which maps
(\ref{ass1}) into
(\ref{wave})) was presented in \cite{Fer96}, see also \cite{Fer97}.

\section{Examples}

Particular exact solutions of the system (\ref{hydro0}) correspond to
special 
surfaces with a flat centroaffine metric.

{\bf Example 1.} Solutions of system (\ref{hydro}) with a degenerate
hodograph
(rarefaction waves) are of the form
$$
p=\varphi (x+cy), ~~~ a=c\varphi (x+cy) +\mu, ~~~
b=c^2\varphi(x+cy) -\frac{1}{\mu}, ~~~ q=c^3\varphi
(x+cy)+c^2\mu-\frac{c}{\mu},
$$
where $c$ and $\mu$ are arbitrary constants and $\varphi$ is an arbitrary
function of $x+cy$. The corresponding potential $f$ is
$$
f(x, y)= \psi (x+cy)+\frac{\mu}{2} x^2y-\frac{1}{2\mu}xy^2 +
\frac{1}{6}(c^2\mu-c/\mu)y^3,
$$
here $\psi '''=\varphi$. The corresponding equations (\ref{paramass}) (where
we have set $\lambda =1$)
take the form
$$
\begin{array}{c}
{\bf r}_{xx}=(c\varphi +\mu) \ {\bf r}_x+\varphi \ {\bf r}_y, \\
\ \\
{\bf r}_{xy}=(c^2\varphi -\frac{1}{\mu}) \ {\bf r}_x+
(c\varphi +\mu) \ {\bf r}_y+ {\bf r}, \\
\ \\
{\bf r}_{yy}=(c^3\varphi +c^2\mu -\frac{c}{\mu}) \ {\bf r}_x+
(c^2\varphi -\frac{1}{\mu}) \ {\bf r}_y.
\end{array}
$$
In the new parametrization $t=(x+cy)/{2}, \ s=(x-cy)/{2}$, these equations
can be rewritten as
\begin{equation}
\begin{array}{c}
{\bf r}_{ss}=(2\mu+\frac{1}{c\mu}) \ {\bf r}_s-\frac{2}{c} \ {\bf r}, \\
\ \\
{\bf r}_{st}=\frac{1}{c\mu} \ {\bf r}_t \\
\ \\
{\bf r}_{tt}=(4c\varphi +2\mu -\frac{2}{c\mu}) \ {\bf r}_t
-\frac{1}{c\mu} \ {\bf r}_s+\frac{2}{c}\ {\bf r},
\end{array}
\label{Ex1}
\end{equation}
implying
$$
({\bf r}-c\mu {\bf r}_s)_s=2\mu ({\bf r}-c\mu {\bf r}_s), ~~~
({\bf r}-c\mu {\bf r}_s)_t=0.
$$
Therefore,
$$
{\bf r}= e^{2\mu s}\ {\bf R_0}+e^{\frac{s}{c\mu}} \ {\bf R}(t)
$$
 where ${\bf R_0}$ is a constant vector and ${\bf R}(t)$ satisfies the ODE
$$
{\bf R} ''=\left(4c\varphi +2\mu -\frac{2}{c\mu}\right) \ {\bf R}'+
\left(\frac{2}{c}-\frac{1}{c^2\mu ^2}\right)\ {\bf R}.
$$
Hence, ${\bf R}(t)$ is a planar curve (which can be arbitrary since
 $\varphi$ is an arbitrary function of $t$), and the intersection of our
surface 
with the plane spanned by ${\bf R_0}$ and ${\bf R}$ has parametric equation
$(e^{2\mu s}, \ e^{\frac{s}{c\mu}})$.

\bigskip

{\bf Example 2.} One can show that, in  appropriate asymptotic
parametrization,
 surfaces of revolution with a flat centroaffine
metric correspond to the solution
$$
a=y\alpha(\xi), ~~~ b=x\beta(\xi), ~~~ p=y^3\rho(\xi), ~~~ q=x^3\gamma(\xi)
$$
where $\xi=xy$ and the functions $\alpha, \beta, \rho, \gamma$ satisfy the
ODEs
$$
\alpha '=3\rho+\xi \rho ', ~~~ \alpha +\xi \alpha '=\beta+\xi \beta ', ~~~
\beta '=3\gamma+\xi \gamma ', ~~~ \xi^3\rho \gamma-\xi \alpha \beta =1.
$$
The general solution of this system can be parametrized in the form
$$
\alpha=2G+\xi G', ~~~ \beta=2H+\xi H', ~~~ \rho=G', ~~~ \gamma = H'
$$
where the functions $G(\xi)$ and $H(\xi)$ solve the quadratic system
$$
\xi^2(H-G)=\epsilon \xi +\nu, ~~~ 2\xi^2HG=\mu-\xi,
$$
$\epsilon, \nu, \mu$ being arbitrary constants. In a different
parametrization,
these surfaces were also discussed in \cite{Scharlach}.
Equivalently, the above solutions can be characterized in terms of the
corresponding 
potential
$$
f(x, y)=F(\xi)+(s_1\xi +s_2) \ln x+(r_1\xi +r_2) \ln y,
$$
where $\xi =xy$ and $s_1, s_2, r_1, r_2$ are arbitrary constants.
Calculation of the 
third order derivatives of $f$ gives
$$
\begin{array}{c}
a=f_{xxy}=y(\xi  F'''+2F''+\frac{s_1}{\xi}), ~~~
b=f_{xyy}=x(\xi  F'''+2F''+\frac{r_1}{\xi}), \\
\ \\
p=f_{xxx}=y^3(F'''-\frac{s_1}{\xi ^2}+2\frac{s_2}{\xi ^3}), ~~~
q=f_{yyy}=x^3(F'''-\frac{r_1}{\xi ^2}+2\frac{r_2}{\xi ^3})
\end{array}
$$
(so that $H=F''+r_1/\xi-r_2/\xi^2   , \ G=F''+s_1/\xi-s_2/\xi^2$), and
the substitution into  (\ref{ass1}) implies the quadratic
equation for $F''$,
$$
\xi^2(F'')^2+(r_1+s_1)\xi F''-(r_2+s_2) F''
-\frac{r_1s_2+r_2s_1}{\xi}+\frac{r_2s_2}{\xi^2}+\frac{1}{2}\xi+s=0,
$$
where $s$ is a constant of integration.

\bigskip 

\noindent
Some further explicit solutions of
equation (\ref{ass1})  can be sought in the form
\begin{equation}
f=y^{3/2} F(x),
\label{ansatz1}
\end{equation}
where $F$ satisfies the ODE $(F^2)'''=-16/3$, implying
$F(x)=\sqrt {-\frac{8}{9} x^3+\alpha x^2+\beta x +\gamma}$,
$\alpha, \beta, \gamma$ being arbitrary constants.
Another possible choice is
\begin{equation}
f=y^3 F(\xi), ~~~ \xi=x/y,
\label{ansatz2}
\end{equation}
implying the following third order  ODE for $F$,
$$
6FF'''-4\xi F'F'''+2\xi (F'')^2-2F'F''=1.
$$ 
A more general self-similar substitution is
$$
f=(xy)^{3/2} F(xy^{\mu})
$$
where $\mu = {\rm const}$. For $\mu=0$ and $\mu=-1$ this reduces to
(\ref{ansatz1}) and (\ref{ansatz2}), respectively.

\section{The loop group formulation of surfaces with a flat centroaffine
metric}

The results of this section are due to A. Bobenko.
Substituting $\psi=e^{\frac{2}{3} \lambda f_{xy}} \ \varphi$ into
(\ref{trace}),
 we obtain the equivalent
traceless representation
\begin{equation}
\varphi_x=\lambda \left(
\begin{array}{ccc}
-\frac{2}{3} a&1&0\\
0&\frac{1}{3} a&p\\
1&b&\frac{1}{3} a
\end{array} \right) \varphi, ~~~~
\varphi_y=\lambda \left(
\begin{array}{ccc}
-\frac{2}{3} b&0&1\\
1&\frac{1}{3} b&a\\
0&q&\frac{1}{3} b
\end{array} \right) \varphi,
\label{traceless}
\end{equation}
where we have set $f_{xxx}=p, \ f_{xxy}=a, \ f_{xyy}=b, \ f_{yyy}=q$.
Let us introduce the loop group of $3\times 3$ matrices
$$
G=\left\{ \varphi (\lambda)~\vert ~~ T(\varphi^t(-\lambda))^{-1} T=\varphi
(\lambda), ~ 
{\rm det} \varphi (\lambda)=1\right\}, ~~~~~
T=\left(
\begin{array}{ccc}
1&0&0\\
0&0&1\\
0&1&0\end{array} \right),
$$
with the corresponding loop algebra
$$
g =\left\{ A(\lambda)~\vert ~~ -TA^t(-\lambda) T=A (\lambda), ~
{\rm tr} A (\lambda)=0\right\}.
$$

{\bf Lemma.} Let $\varphi (t, \lambda) \to G$ be such that
$\varphi _t \varphi ^{-1}$ is linear in $\lambda$, and
$\varphi _t \varphi ^{-1}\vert _{\lambda =0} =0$. Then
$$
\varphi _t \varphi ^{-1} (t, \lambda)=\lambda
\left(
\begin{array}{ccc}
-2A&C&D\\
D&A&P\\
C&B&A
\end{array} \right).
$$
This readily implies the following

{\bf Proposition.} Let $\varphi (u, v, \lambda) \to G$ be an immersion
of an open domain of $R^2$ such that
$d\varphi \varphi ^{-1}$ is linear in $\lambda$, and
$d\varphi \varphi ^{-1}\vert _{\lambda =0} =0$. Then there exists a
change of variables $(u, v) \to (x, y)$ bringing the frame equations
into the form (\ref{traceless}).

\section{Projective generalizations of hypersurfaces with a flat
centroaffine metric}

There exist two particularly interesting projectively invariant
classes of hypersurfaces containing projective transforms of hypersurfaces
with
flat centroaffine metric.

\bigskip 
\noindent {\bf Hypersurfaces with conformally flat second fundamental form}
were investigated in \cite {Ak1}, \cite{Ak2} and \cite{Kon} (notice that
this 
property makes sense only for hypersurface of dimension greater than two).
Among 
 previosly known examples one should mention the second order
envelopes of one-parameter
families of nondegenerate hyperquadrics
(in particular, projective transforms of surfaces of revolution)
as well as  special projections of the Segre variety $P^2 \times P^2 \subset
P^8$ 
into $P^5$. We point out that projective transforms of hypersurfaces with
flat 
centroaffine metric provide further examples of this type. Indeed, the
conformal class of the second fundamental form of a hypersurface with flat
centroaffine metric
is $\eta_{ij}du^idu^j$, which is manifestly conformally flat.

\bigskip
\noindent {\bf Surfaces in $P^3$ possessing asymptotic deformations which
induce rescalings of the projective metric and the Darboux cubic form}.
Let us recall that the position vector ${\bf R}$ of a surface $M^2$
in the projective space $P^3$ satisfies the linear system
\begin{equation}
\begin{array}{c}
{\bf R}_{xx}=\beta \ {\bf R}_y+\frac{1}{2}(V-\beta_y) \ {\bf R} \\
\ \\
{\bf R}_{yy}=\gamma \ {\bf R}_x+\frac{1}{2}(W-\gamma_x) \ {\bf R}
\end{array}
\label{proj1}
\end{equation}
where $\beta, \gamma, V, W$ are functions of the asymptotic coordinates
$x$ and  $y$ satisfying the compatibility
conditions \cite[p.\ 120]{Lane}
\begin{equation}
\begin{array}{c}
\beta_{yyy}-2\beta_yW-\beta W_y=
\gamma_{xxx}-2\gamma_xV-\gamma V_x\\
\ \\
W_x=2\gamma \beta_y+\beta \gamma_y \\
\ \\
V_y=2\beta \gamma_x+\gamma \beta_x.
\end{array}
\label{proj2}
\end{equation}
The most important projective invariants are the projective metric
$$
2 \beta \gamma\,dxdy
$$
and the conformal class of the Darboux cubic form
$$
\beta \,dx^3+\gamma \,dy^3.
$$
Paticularly interesting classes of projective surfaces are characterized by
additional relations among  $\beta, \gamma, V$ and $ W$, for instance

\noindent -- isothermally-asymptotic surfaces ($\beta=\gamma$);

\noindent -- projectively applicable surfaces ($\beta_y=\gamma_x$);

\noindent -- surfaces of Jonas ($\beta_x=\gamma_y$);

\noindent etc (see \cite{Bol}, \cite{Lane}, \cite{Fer00}).

\noindent Here we consider another class of projective surfaces
characterised by
the additional relation
\begin{equation}
\beta_{yyy}=\gamma_{xxx},
\label{3}
\end{equation}
which allows one to introduce the spectral parameter $\lambda$ into
(\ref{proj1}):
\begin{equation}
\begin{array}{c}
{\bf R}_{xx}=\lambda \beta \ {\bf R}_y+\frac{1}{2}(\lambda^2 V-
\lambda \beta_y) \ {\bf R} \\
\ \\
{\bf R}_{yy}=\lambda \gamma \ {\bf R}_x+\frac{1}{2}(\lambda^2 W-
\lambda \gamma_x) \ {\bf R}.
\end{array}
\label{proj3}
\end{equation}
Since the projective metric and the conformal class of the cubic form of the
one-parameter family
of surfaces $M^2_{\lambda}$ defined by (\ref{proj3}) are
$$
2 \lambda^2 \beta \gamma\,dxdy
$$
and 
$$
\beta \,dx^3+\gamma \,dy^3,
$$
respectively, one can say that the surfaces satisfying (\ref{3}) possess
one-parameter families of asymptotic deformations which rescale the
projective metric and preserve the conformal class of the
Darboux cubic form.
As $\lambda$ varies, we continuosly change from the initial surface
($\lambda=1$)
to a quadric ($\lambda=0$), and then to the dual surface ($\lambda=-1$).
Surfaces satisfying the condition
(\ref{3}) were discussed in \cite{Cech} and \cite{Kaucky}.
Changing in the equations (\ref{Lax1}) from ${\bf r}$ to
${\bf R}= e^{-\frac{1}{2}\lambda f_{xy}}{\bf r}$,  we readily rewrite
$(\ref{Lax1})_1$ and $(\ref{Lax1})_3$ in the form (\ref{proj3}), where
$$
\beta=f_{xxx}, ~~~ \gamma=f_{yyy}, ~~~
V=\frac{1}{2}f_{xxy}^2+f_{xxx}f_{xyy}, ~~~
W=\frac{1}{2}f_{xyy}^2+f_{yyy}f_{xxy}
$$
and $f$ satisfies the associativity equation (\ref{ass1}). Therefore,
projective transforms of surfaces (\ref{Lax1}) form a proper subclass of
the integrable class of surfaces characterized by the condition (\ref{3}).

\section{Acknowledgements}

It is a great pleasure to thank A.~I.  Bobenko for numerous discussions and
for proving the material  presented in  sect. 7. This research was partially
supported by   SFB 288  and the EPSRC grant GR/N30941.


\begin{thebibliography}{99}
\addcontentsline{toc}{section}{References}

\bibitem{Ak1} Akivis M.A. and Goldberg V.V.,
Conformal differential geometry and its
generalizations, John Wiley and Sons, 1996.

\bibitem{Ak2} Akivis M.A. and Konnov V.V., Local aspects in conformal
structure theory, 
Russian Math. Surveys {\bf 48}, no. 1 (1993) 1-35.


\bibitem{Bol} Bol G., Projektive Differentialgeometrie,
  G\"ottingen, 1954.


\bibitem{Cartan} Cartan~E., Sur la d\'eformation projective des surfaces,
Ann. de l'Ecole Normale Sup., (3) 37 (1920) 259-356.

\bibitem{Cech} \^Cech E., Sur les correspondances asymptotiques entre
deux surfaces, Rendiconti  Accademia Nazionale dei Lincei, Roma {\bf 6} ser.
8
(1928) 484-486; 552-554.

\bibitem{Dubrovin} Dubrovin B.~A., Geometry of 2D topological field
theories, Lecture Notes in Mathematics, V.1620, Berlin,
Springer, 120-348.

\bibitem{Dryuma} Dryuma V., On the Riemannian and Einstein-Weyl geometry in
the theory 
of second order ordinary differential equations, arXiv:math.DG/0104278.

\bibitem{Fer96} Ferapontov E.~V. and Mokhov O.~I., Equations of
associativity
of two-dimensional topological field theory as integrable Hamiltonian
nondiagonalisable systems of hydrodynamic type, Funkt. Anal. and it's Appl.,
{\bf 30}, N3 (1996) 62-72.



\bibitem{Fer99}  Ferapontov E.~V. and Schief W.~K., Surfaces of Demoulin:
  Differential geometry, B\"acklund transformation and Integrability,
  Journal of Geometry and Physics {\bf 30} (1999) 343-363.

\bibitem{Fer00}  Ferapontov~E.V. Integrable systems in projective
differential 
geometry, Kyushu J. Math. {\bf 54} (2000) 183-215.

\bibitem{Fer97} Ferapontov E.~V.,  Galvao C.~A.~P.,  Mokhov O.~I. and Nutku
Y. ,
  Bi-Hamiltonian structure of equations of
  associativity in 2-d topological field theory, Comm. in Math. Phys., {\bf
186} (1997)
  649-669.

\bibitem{Finikov} Finikov S.~P., Projective Differential Geometry,
Moscow-Leningrad, 1937.


\bibitem{Fubini} Fubini G. and \^Cech E., Geometria Proiettiva
Differenziale,
Bologna: Zanichelli, 1926.



\bibitem{Kaucky} Kauck\'y J., Sur les transformations asymptotiques
d'une surface non d\'evelopable en elle
m\^eme dans l'espace projectif $S_3$, M\'em. Boh\^ eme {\bf 16}
(1932) 1-35.

\bibitem{Kon}  Konnov V.V., Asymptotic pseudoconformal structure on a
four-dimensional 
hypersurface and its completely isotropic two-dimensional submanifolds,
Russian Math. {\bf 36}, no. 6 (1992) 67-74.


\bibitem{Lane} Lane E., Treatise on Projective Differential Geometry,
The Univ. of Chicago Press, 1942.

\bibitem{Wang} Liu H. and Wang C., The centroaffine Chebyshev operator,
Results in Math. {\bf 27} (1995) 77-92.

\bibitem{Loewner} Loewner C. and Nirenberg L., Partial differential
equations 
invariant under conformal or projective transformations, in Contributions to
Analysis, Academic Press (1974) 245-272.

\bibitem{Loftin} Loftin J.~C., Riemannian metrics on locally projectively
flat
manifolds, arXiv:math.DG/0108218

\bibitem{Mokhov} Mokhov O.I., Symplectic and Poisson geometry on loop spaces
of 
manifolds and nonlinear
equations, AMS Translations ser. 2 {\bf 170} (1995).


\bibitem{Sasaki1} Sasaki T., Projective differential geometry and
linear homogeneous differential equations, Rokko lectures in Math.,
Kobe University, Japan, {\bf 5} (1999) 1-115, ISBN 4-907719-05-1.

\bibitem{Scharlach} Scharlach~C., Some results in centroaffine differential
geometry,
in: Geometry and Topology of Submanifolds, {\bf 4}, ed. F. Dillen and
L. Verstraelen, World Scientific, Singapore, (1992) 198-206.

\bibitem{Shirokov} Shirokov P.~A. and Shirokov A.~P., Affine Differential
geometry,
Moscow, 1959.

\bibitem{Strachan} Strachan~I.A.B., Frobenius submanifolds,
J. Geom. Phys. {\bf 38}  (2001) 285--307.

\bibitem{Wilczynski} Wilczynski E.~I., Projective-differential geometry of
curved surfaces, Trans. AMS {\bf 8} (1907) 233-260; {\bf 9}
 (1908) 79-120, 293-315.

\bibitem{Yoshida} Yoshida M., Fuchsian Differential Equations, Vieweg, 1987.

\end{thebibliography}
\end{document}